\documentclass{article}
\usepackage{graphicx}
\usepackage{amsmath}
\usepackage{amsfonts}
\usepackage{amssymb}

\begin{document}
\LaTeX{}{}

\bigskip

\begin{center}
\ \ \textbf{Score lists in [h-k]-bipartite hypertournaments}
\end{center}

$\bigskip$

\begin{center}
S. Pirzada$^{1}$, T. A. Chishti$^{2}$, T. A. Naikoo$^{3}$
\end{center}

\bigskip

$^{1,3}$Department of Mathematics, University of Kashmir, India

$^{1}$Email: sdpirzada@yahoo.co.in \ \ \ \ \ \ \ \ \ 

$^{3}$Email: tariqnaikoo@rediffmail.com

$^{2}$Centre of Distance Education, University of Kashmir, India\strut

\bigskip$^{2}$Email: chishtita@yahoo.co.in\ 


\bigskip

\textbf{ABSTRACT. }Given non-negative integers m, n, h and k with $\ m\geq
h>1\ $ and $\ n\geq k>1,\ $an [h-k]-bipartite\ hypertournament on $\ m+n$
\ vertices is a triple $\left(  U,V,A\right)  $, where U and V are two sets of
vertices with $\left|  U\right|  =m$ and $\ \left|  V\right|  =n,$ and $A$ is
a set of \ $\left(  h+k\right)  -$ tuples of vertices, called arcs, with
exactly $h$ vertices from $U$ \ and exactly $k$ vertices from $V$ ,such that
any $h+k$  subsets $\ U_{1}\cup V_{1}$ of \ $U\cup V,\ A$ contains exactly one
of the $\left(  h+k\right)  !\ \ \left(  h+k\right)  -$tuples whose entries
belong to $U_{1}\cup V_{1}.$ We obtain necessary and sufficient conditions for
a pair of non-decreasing sequences of non-negative integers to be the losing
score lists or score lists of some$[h-k]-$bipartite \ hypertournament.\bigskip

\begin{center}
\textbf{1. Introduction\bigskip}
\end{center}

Hypergraphs are generalization of graphs [3]. While edges of a graph are pairs
of vertices of the graph, edges of a hypergraph are subsets of the vertex set,
consisting of at least two vertices. An edge consisting of k vertices is
called a k-edge. A k-hypergraph is a hypergraph all of whose edges are
k-edges. A k-hypertournament is a complete k-hypergraph with each k-edge
endowed with an orientation, that is, a linear arrangement of the vertices
contained in the hyperedge.

Instead of scores of vertices in a tournament, Zhou et al. [8] considered
scores and losing scores of vertices in a k-hypertournament, and derived a
result analogous to Landau's theorem [6]. The score s(v$_{i}$) or s$_{i}$ of a
vertex v$_{i}$ is the number of arcs containing v$_{i}$\ and in which v$_{i}%
$\ is not the last element, and the losing score r(v$_{i}$) or r$_{i}$\ of a
vertex v$_{i}$\ is the number of arcs containing v$_{i}$\ and in which v$_{i}%
$\ is the last element. The score sequence (losing score sequence) is formed
by listing the scores (losing scores) in non-decreasing order.

We note that for two integers p and q,

\begin{center}
$\left(
\begin{array}
[c]{c}%
p\\
q
\end{array}
\right)  $ $=\left\{
\begin{array}
[c]{c}%
\frac{p!}{q!(p-q)!},\qquad p\geq q,\\
0,\qquad\ \ \ \ \ \ \ \ \ p<q.
\end{array}
\right.  $
\end{center}

The following characterizations of score sequences and losing score sequences
in k-hypertournaments can be found in Zhou et al. [8].\bigskip

\textbf{Theorem 1.1.} Given two non-negative integers n and k with $n\geq
k>1,$\ a non-decreasing sequence R = [r$_{1}$, r$_{2}$ ,..., r$_{n}$] of
non-negative integers is a losing score sequence of some k-hypertournament if
and only if for each j,

\begin{center}
$\sum_{i=1}^{j}r_{i}\geq\left(
\begin{array}
[c]{c}%
j\\
k
\end{array}
\right)  ,$\ 
\end{center}

with equality when j = n.\bigskip

\textbf{Theorem 1.2.} Given non-negative integers n and k with $n\geq k>1,$\ a
non-decreasing sequence S = [s$_{1}$, s$_{2}$ ,..., s$_{n}$] of non-negative
integers is a score sequence of some k-hypertournament if and only if for each j,

\begin{center}
$\sum_{i=1}^{j}s_{i}\geq j\left(
\begin{array}
[c]{c}%
n-1\\
k-1
\end{array}
\right)  +\left(
\begin{array}
[c]{c}%
n-j\\
k
\end{array}
\right)  -\left(
\begin{array}
[c]{c}%
n\\
k
\end{array}
\right)  ,$
\end{center}

with equality when j = n.\ \bigskip\ \ 

Bang and Sharp [1] proved Landau's theorem using Hall's theorem on a system of
distinct representatives of a collection of sets. Based on Bang and Sharp's
ideas, Koh and Ree [5] have given a different proof of Theorems 1.1 and
1.2.\bigskip\ Some more results on scores of k-hypertournaments can be found
in [4, 7].

Bipartite hypergraphs are generalization of bipartite graphs. If U = \{u$_{1}
$, u$_{2}$ ,..., u$_{m}$\} and V = \{v$_{1}$, v$_{2}$ ,..., v$_{n}$\} are
vertex sets, then the edge of a bipartite hypergraph is a subset of the vertex
sets, containing at least one vertex from U and at least one vertex from V. If
an edge has exactly h vertices from U and exactly k vertices from V, it is
called an [h-k]-edge. An [h-k]-bipartite hypergraph is a bipartite hypergraph
all of whose edges are [h-k]-edges. An [h-k]-bipartite hypertournament is a
complete [h-k]-bipartite hypergraph with each [h-k]-edge endowed with an
orientaion, that is, a linear arrangement of the vertices contained in the hyperedge.

Equivalently, given non-negative integers m, n, h and k with $\ m\geq h>1\ $
and $\ n\geq k>1,\ $an [h-k]-bipartite hypertournament of order m x n consists
of two vertex sets U and V with $\left|  U\right|  =m$ and $\ \left|
V\right|  =n,$ together with an arc set E, a set of\ $\left(  h+k\right)  $
tuples of vertices, with exactly h vertices from U and exactly k vertices from
V, called arcs, such that for any h-subset U$_{1}$ of U and k-subset V$_{1}$
of V, E contains exactly one of the (h+k)! \ (h+k)-tuples whose h entries
belong to U$_{1}$ and k entries belong to V$_{1}$. Let e = (u$_{1}$, u$_{2}%
$,..., u$_{h}$, v$_{1}$, v$_{2}$,..., v$_{k}$) be an arc in H and i
$<$%
j, we denote e(u$_{i}$, u$_{j}$) = (u$_{1}$,..., u$_{j}$,..., u$_{i}$, v$_{1}%
$,..., v$_{k}$), that is, the new arc obtained from e by interchanging u$_{i}$
and u$_{j}$ in e. Similarly, we can have new arcs of the form e(v$_{i}$,
v$_{j}$) and e(u$_{i}$, v$_{j}$).

For a given vertex u$_{i}$\ $\in$\ U, the score $d_{H}^{+}(u_{i})$\ (or simply
$d^{+}(u_{i})$) is the number of [h-k]-arcs containing u$_{i}$\ and in which
u$_{i}$\ is not the last element. The losing score $d_{H}^{-}(u_{i}) $\ (or
simply $d^{-}(u_{i})$)\ is the number of [h-k]-arcs containing u$_{i}$\ and in
which u$_{i}$\ is the last element. Similarly, we define by $d_{H}^{+}(v_{j}%
)$\ and $d_{H}^{-}(v_{j})$\ respectively as the score and losing score of a
vertex v$_{j}$ $\in$\ V. The losing score lists of an [h-k]-bipartite
hypertournament is a pair of non-decreasing sequences of non-negative integers
A = [a$_{1}$, a$_{2}$ ,..., a$_{m}$] and B = [b$_{1} $, b$_{2}$ ,..., b$_{n}%
$], where a$_{i}$ is a losing score of some vertex u$_{i} $ $\in$\ U and
b$_{j}$ is a losing score of some vertex v$_{j} $ $\in$\ V. Similarly, the
score lists are formed by listing the scores in non-decreasing order, and we
denote these by C = [c$_{1}$, c$_{2}$ ,..., c$_{m}$] and D = [d$_{1}$, d$_{2}$
,..., d$_{n}$].\bigskip

\begin{center}
\textbf{2. Main results} \bigskip
\end{center}

The following two Theorems are the main results and provide a characterization
of losing score lists and score lists in [h-k]-bipartite hypertournaments.\bigskip

\textbf{Theorem 2.1. }Given non-negative integers m, n, h and k with $\ m\geq
h>1\ $ and $\ n\geq k>1,$the non-decreasing sequences A = [a$_{i}$]$_{1}^{m}$
and B = [b$_{j}$]$_{1}^{n}$ of non-negative integers are the losing score
lists of an [h-k]-bipartite hypertournament if and only if for each p and q,

\begin{center}
$\ \ \ \ \ \ \sum_{i=1}^{p}a_{i}+\sum_{j=1}^{q}b_{j}\geq\left(
\begin{array}
[c]{c}%
p\\
h
\end{array}
\right)  \left(
\begin{array}
[c]{c}%
q\\
k
\end{array}
\right)  ,\qquad\ \ \ \ \ \ \ \ \ \ \ \ \ \ \ \ \ \ \ \ (1)$
\end{center}

with equality when p = m and q = n.\bigskip

\textbf{Theorem 2.2. }Given non-negative integers m, n, h and k with $\ m\geq
h>1\ $ and $\ n\geq k>1,$the non-decreasing sequences C = [c$_{i}$]$_{1}^{m}$
and D = [d$_{j}$]$_{1}^{n}$ of non-negative integers are the score lists of an
[h-k]-bipartite hypertournament if and only if for each p and q,

$\ \sum_{i=1}^{p}c_{i}+\sum_{j=1}^{q}d_{j}\geq p\left(
\begin{array}
[c]{c}%
m-1\\
h-1
\end{array}
\right)  \left(
\begin{array}
[c]{c}%
n\\
k
\end{array}
\right)  +q\left(
\begin{array}
[c]{c}%
m\\
h
\end{array}
\right)  \left(
\begin{array}
[c]{c}%
n-1\\
k-1
\end{array}
\right)  +$

$\smallskip
\ \ \ \ \ \ \ \ \ \ \ \ \ \ \ \ \ \ \ \ \ \ \ \ \ \ \ \ \ \ \ \ \ \ \ \ \ \ \left(
\begin{array}
[c]{c}%
m-p\\
h
\end{array}
\right)  \left(
\begin{array}
[c]{c}%
n-q\\
k
\end{array}
\right)  -\left(
\begin{array}
[c]{c}%
m\\
h
\end{array}
\right)  \left(
\begin{array}
[c]{c}%
n\\
k
\end{array}
\right)  ,\ \ \ \ \ \ \ \ \ \ \ \ \ \ \ \ \ \ \ \ (2)$

with equality when p = m and q = n.\bigskip\bigskip

In order to prove Theorem 2.1 and Theorem 2.2, we require the following
Lemmas. We note that in an [h-k]-bipartite hypertournament H there are exactly
$\left(
\begin{array}
[c]{c}%
m\\
h
\end{array}
\right)  \left(
\begin{array}
[c]{c}%
n\\
k
\end{array}
\right)  $ arcs, and in each arc, only one vertex is at the last entry. Therefore,

\begin{center}
$\ \sum_{i=1}^{m}d_{H}^{-}(u_{i})+\sum_{j=1}^{n}d_{H}^{-}(v_{j})=\left(
\begin{array}
[c]{c}%
m\\
h
\end{array}
\right)  \left(
\begin{array}
[c]{c}%
n\\
k
\end{array}
\right)  .\bigskip$
\end{center}

\textbf{Lemma 2.1. }If H is an [h-k]-bipartite hypertournament of order m x n
with score lists A = [c$_{i}$]$_{1}^{m}$ and B = [d$_{j}$]$_{1}^{n}$, then

\begin{center}
$\sum_{i=1}^{m}c_{i}+\sum_{j=1}^{n}d_{j}=(h+k-1)\left(
\begin{array}
[c]{c}%
m\\
h
\end{array}
\right)  \left(
\begin{array}
[c]{c}%
n\\
k
\end{array}
\right)  .$
\end{center}

\textbf{Proof. }Obviously, m $\geq$\ h and n $\geq$\ k. If a$_{i}$ is the
losing score of u$_{i}$ $\in$\ U and b$_{j}$ is the losing score of v$_{j}$
$\in$\ V, then

\begin{center}
$\ \ \ \sum_{i=1}^{m}a_{i}+\sum_{j=1}^{n}b_{j}\geq\left(
\begin{array}
[c]{c}%
m\\
h
\end{array}
\right)  \left(
\begin{array}
[c]{c}%
n\\
k
\end{array}
\right)  .$
\end{center}

Now, there are $\left(
\begin{array}
[c]{c}%
m-1\\
h-1
\end{array}
\right)  \left(
\begin{array}
[c]{c}%
n\\
k
\end{array}
\right)  $\ arcs containing a vertex u$_{i}$ $\in$\ U, and $\left(
\begin{array}
[c]{c}%
m\\
h
\end{array}
\right)  \left(
\begin{array}
[c]{c}%
n-1\\
k-1
\end{array}
\right)  $\ arcs containing a vertex v$_{j}$ $\in$\ V. Therefore,

$\sum_{i=1}^{m}c_{i}+\sum_{j=1}^{n}d_{j}=\sum_{i=1}^{m}\left(
\begin{array}
[c]{c}%
m-1\\
h-1
\end{array}
\right)  \left(
\begin{array}
[c]{c}%
n\\
k
\end{array}
\right)  +$

$\ \ \ \ \ \ \ \ \ \ \ \ \ \ \ \ \ \ \ \ \ \ \ \ \ \ \ \ \sum_{j=1}%
^{n}\left(
\begin{array}
[c]{c}%
m\\
h
\end{array}
\right)  \left(
\begin{array}
[c]{c}%
n-1\\
k-1
\end{array}
\right)  -\left(
\begin{array}
[c]{c}%
m\\
h
\end{array}
\right)  \left(
\begin{array}
[c]{c}%
n\\
k
\end{array}
\right)  $

$\smallskip\qquad\ \ \ \ \ \ \ =m\left(
\begin{array}
[c]{c}%
m-1\\
h-1
\end{array}
\right)  \left(
\begin{array}
[c]{c}%
n\\
k
\end{array}
\right)  +n\left(
\begin{array}
[c]{c}%
m\\
h
\end{array}
\right)  \left(
\begin{array}
[c]{c}%
n-1\\
k-1
\end{array}
\right)  -\left(
\begin{array}
[c]{c}%
m\\
h
\end{array}
\right)  \left(
\begin{array}
[c]{c}%
n\\
k
\end{array}
\right)  \smallskip$

$\ \ \ \ \ \ \ \ \ \ \ \ \ \ =(h+k-1)\left(
\begin{array}
[c]{c}%
m\\
h
\end{array}
\right)  \left(
\begin{array}
[c]{c}%
n\\
k
\end{array}
\right)  .\qquad\ \ \ \ \
\begin{array}
[c]{c}
\end{array}
\bigskip$

\textbf{Lemma 2.2. }If A = [a$_{1}$, a$_{2}$ ,..., a$_{m}$] and B = [b$_{1}$,
b$_{2}$ ,..., b$_{n}$] are losing score lists of an [h-k]-bipartite
hypertournament H, and if a$_{i}$
$<$%
a$_{j}$, then A$^{/}$ = [a$_{1}$, a$_{2}$ ,..., a$_{i}$+1,..., a$_{j}%
$-1,...,a$_{m}$] and B are losing score lists of some [h-k]-bipartite hypertournament.

\textbf{Proof. }Let A and B be the\ losing score lists of an [h-k]-bipartite
hypertournament H with vertex sets U = \{u$_{1}$, u$_{2}$,..., u$_{m}$\} and V
= \{v$_{1}$, v$_{2}$,..., v$_{n}$\} so that $d^{-}(u_{i})=a_{i}$\ and
$d^{-}(v_{j})=b_{i}$ \ $(1\leq i\leq m,~1\leq j\leq n)$.

If there is an [h-k]-arc e containing both u$_{i}$ and u$_{j}$ with u$_{j}$ as
the last element in e, let e$^{/}$ = (u$_{i}$, u$_{j}$) and H$^{/}$ = (H-e)
$\cup$\ e$^{/}$. Clearly A$^{/}$ and B are the losing score lists of H$^{/}$.

Now, assume that for every arc e containing both u$_{i}$ and u$_{j}$,
\ u$_{j}$ is not the last element in e. Since a$_{i}$
$<$%
a$_{j}$, there exist two [h-k]-arcs e$_{1}$ = (w$_{1}$, w$_{2}$,..., w$_{l-1}%
$, u$_{i} $, w$_{l}$,..., w$_{h-1}$, z$_{1}$, z$_{2}$,..., z$_{k}$) and
e$_{2}$ = (w$_{1}^{/}$, w$_{2}^{/}$,..., w$_{h-1}^{/}$, z$_{1}^{/}$,
z$_{2}^{/}$,..., z$_{k}^{/}$, u$_{j}$) where w's $\in$\ U, z's $\in$\ V,
u$_{i}$ $\notin$\ \{w$_{1}$, w$_{2}$,..., w$_{h-1}$\}, u$_{j}$ $\notin
$\ \{w$_{1}$, w$_{2}$,..., w$_{h-1}$\} and (w$_{1}^{/}$, w$_{2}^{/}$,...,
w$_{h-1}^{/}$, z$_{1}^{/}$, z$_{2}^{/}$,..., z$_{k}^{/}$) is a permutation of
(w$_{1}$, w$_{2} $,..., w$_{h-1}$, z$_{1}$, z$_{2}$,..., z$_{k}$).

Now, let e$_{1}^{/}$ = e$_{1}$(u$_{i}$, x) and e$_{2}^{/}$ = e$_{2}$(u$_{j}$,
y) where x is any one from \{w$_{1}$, w$_{2}$,..., w$_{h-1}$, z$_{1}$, z$_{2}%
$,..., z$_{k}$\}and y is any one from \{w$_{1}^{/}$, w$_{2}^{/}$,...,
w$_{h-1}^{/}$, z$_{1}^{/}$, z$_{2}^{/}$,..., z$_{k}^{/}$\}. Take H$^{/}$ =
(H-(e$_{1}$ $\cup$\ e$_{2}$)) $\cup$\ (e$_{1}^{/}$ $\cup$\ e$_{2}^{/}$). Then,
A$^{/}$ and B are the score lists of H$^{/}$.$%
\begin{array}
[c]{c}
\end{array}
\bigskip$

\textbf{Lemma 2.3. }Let A = [a$_{1}$, a$_{2}$ ,..., a$_{m}$] and B = [b$_{1}$,
b$_{2}$ ,..., b$_{n}$] be non-decreasing sequences of non-negative integers
satisfying (1). If $a_{m}<\left(
\begin{array}
[c]{c}%
m-1\\
h-1
\end{array}
\right)  \smallskip\left(
\begin{array}
[c]{c}%
n\\
k
\end{array}
\right)  $, then there exists r (1 $\leq$\ r $\leq$\ m-1) such that A$^{/}$ =
[a$_{1}$, a$_{2}$ ,..., a$_{r}$-1,..., a$_{m}$+1] is non-decreasing and
A$^{/}$ and B satisfy (1).

\textbf{Proof. }Let r be the maximum integer such that \ a$_{r-1}$
$<$%
a$_{r}$ = a$_{r+1}$ = ...= a$_{m-1}$with a$_{0}$ = 0 if r = 1.

To show that A$^{/}$ and B satisfy (1), we need to prove that for each p (r
$\leq$\ p $\leq$\ m-1),\ \ \ \ 

\begin{center}
$\ \ \ \ \ \ \ \ \ \ \ \ \ \ \ \ \ \ \ \ \ \ \ \ \ \ \ \ \ \ \ \ \ \ \ \ \ \ \ \ \ \ \ \ \ \ \ \sum
_{i=1}^{p}a_{i}+\sum_{j=1}^{q}b_{j}>\left(
\begin{array}
[c]{c}%
p\\
h
\end{array}
\right)  \left(
\begin{array}
[c]{c}%
q\\
k
\end{array}
\right)
,~\ \ \ \ \ \ \ \ \ \ \ \ \ \ \ \ \ \ \ \ \ \ \ \ \ \ \ \ \ \ \ \ \ \ \ \ \ \ \ \ \ \ \ \ \ \ \ \ \ \ \ \ \ (3)$%
\ \ \ \ \ 
\end{center}

As \ $a_{m}<\left(
\begin{array}
[c]{c}%
m-1\\
h-1
\end{array}
\right)  \left(
\begin{array}
[c]{c}%
n\\
k
\end{array}
\right)  $, we have\ \ \ \ 

$\ \ \ \sum_{i=1}^{m-1}a_{i}+\sum_{j=1}^{n}b_{j}=\left(
\begin{array}
[c]{c}%
m\\
h
\end{array}
\right)  \left(
\begin{array}
[c]{c}%
n\\
k
\end{array}
\right)  -a_{m}$

$\smallskip
\ \ \ \ \ \ \ \ \ \ \ \ \ \ \ \ \ \ \ \ \ \ \ \ \ \ \ \ \ \ \ >\left(
\begin{array}
[c]{c}%
m\\
h
\end{array}
\right)  \left(
\begin{array}
[c]{c}%
n\\
k
\end{array}
\right)  -\left(
\begin{array}
[c]{c}%
m-1\\
h-1
\end{array}
\right)  \left(
\begin{array}
[c]{c}%
n\\
k
\end{array}
\right)  \smallskip$

$\smallskip
\ \ \ \ \ \ \ \ \ \ \ \ \ \ \ \ \ \ \ \ \ \ \ \ \ \ \ \ \ \ \ =\left[
\left(
\begin{array}
[c]{c}%
m\\
h
\end{array}
\right)  -\left(
\begin{array}
[c]{c}%
m-1\\
h-1
\end{array}
\right)  \right]  \left(
\begin{array}
[c]{c}%
n\\
k
\end{array}
\right)  \smallskip$\ 

$\smallskip
\ \ \ \ \ \ \ \ \ \ \ \ \ \ \ \ \ \ \ \ \ \ \ \ \ \ \ \ \ \ \ =\left(
\begin{array}
[c]{c}%
m-1\\
h
\end{array}
\right)  \left(
\begin{array}
[c]{c}%
n\\
k
\end{array}
\right)  ,~\ $\ \ \ 

This shows that for r = m-1, (3) is true.

Now, assume that r $\leq$\ m-2. Then (3) holds for p = m-1.

If there exists p$_{0}$ (r $\leq\ $p$_{0}$ $\leq$\ m-2) such that

\begin{center}
\bigskip$\sum_{i=1}^{p_{0}}a_{i}+\sum_{j=1}^{q}b_{j}=\left(
\begin{array}
[c]{c}%
p_{0}\\
h
\end{array}
\right)  \left(
\begin{array}
[c]{c}%
q\\
k
\end{array}
\right)  ,$
\end{center}

choose p$_{0}$ as large as possible.

Since $\ \ \ \ \ \ \ \sum_{i=1}^{p_{0}+1}a_{i}+\sum_{j=1}^{q}b_{j}>\left(
\begin{array}
[c]{c}%
p_{0}+1\\
h
\end{array}
\right)  \left(
\begin{array}
[c]{c}%
q\\
k
\end{array}
\right)  ,$

therefore \ \ \ $a_{p_{0}}=a_{p_{0}+1}=\left(  \sum_{i=1}^{p_{0}+1}a_{i}%
+\sum_{j=1}^{q}b_{j}\right)  -\left(  \sum_{i=1}^{p_{0}}a_{i}+\sum_{j=1}%
^{q}b_{j}\right)  $

$\smallskip\ \ \ \ \ \ \ \ \ \ \ \ \ \ \ \ \ \ \ \ \ >\left(
\begin{array}
[c]{c}%
p_{0}+1\\
h
\end{array}
\right)  \left(
\begin{array}
[c]{c}%
q\\
k
\end{array}
\right)  -\left(
\begin{array}
[c]{c}%
p_{0}\\
h
\end{array}
\right)  \left(
\begin{array}
[c]{c}%
q\\
k
\end{array}
\right)  =$\ $\left(
\begin{array}
[c]{c}%
p_{0}\\
h-1
\end{array}
\right)  \left(
\begin{array}
[c]{c}%
q\\
k
\end{array}
\right)  .$

Thus, it follows that$\qquad$

$\sum_{i=1}^{p_{0}-1}a_{i}+\sum_{j=1}^{q}b_{j}=\sum_{i=1}^{p_{0}}a_{i}%
+\sum_{j=1}^{q}b_{j}-a_{p_{0}}$

$\smallskip
\ \ \ \ \ \ \ \ \ \ \ \ \ \ \ \ \ \ \ \ \ \ \ \ \ \ \ \ \ \ \ \ <\left(
\begin{array}
[c]{c}%
p_{0}\\
h
\end{array}
\right)  \left(
\begin{array}
[c]{c}%
q\\
k
\end{array}
\right)  -\left(
\begin{array}
[c]{c}%
p_{0}\\
h-1
\end{array}
\right)  \left(
\begin{array}
[c]{c}%
q\\
k
\end{array}
\right)  $

$\smallskip
\ \ \ \ \ \ \ \ \ \ \ \ \ \ \ \ \ \ \ \ \ \ \ \ \ \ \ \ \ \ \ <\left[
\left(
\begin{array}
[c]{c}%
p_{0}-1\\
h
\end{array}
\right)  +\left(
\begin{array}
[c]{c}%
p_{0}-1\\
h-1
\end{array}
\right)  \right]  \left(
\begin{array}
[c]{c}%
q\\
k
\end{array}
\right)  -\left(
\begin{array}
[c]{c}%
p_{0}\\
h-1
\end{array}
\right)  \left(
\begin{array}
[c]{c}%
q\\
k
\end{array}
\right)  $

$\smallskip
\ \ \ \ \ \ \ \ \ \ \ \ \ \ \ \ \ \ \ \ \ \ \ \ \ \ \ \ \ \ \ <\left[
\left(
\begin{array}
[c]{c}%
p_{0}-1\\
h
\end{array}
\right)  -\left(
\begin{array}
[c]{c}%
p_{0}-1\\
h-2
\end{array}
\right)  \right]  \left(
\begin{array}
[c]{c}%
q\\
k
\end{array}
\right)  $

$\smallskip\ \ \ \ \ \ \ \ \ \ \ \ \ \ \ \ \ \ \ \ \ \ \ \ \ \ \ \ \smallskip
\ \ <\left(
\begin{array}
[c]{c}%
p_{0}-1\\
h
\end{array}
\right)  \left(
\begin{array}
[c]{c}%
q\\
k
\end{array}
\right)  ,$

a contradiction with the hypothesis on A and B. Hence (3) holds.\bigskip
\ \ \ \ \ \ \ \ \ \ \ \ \ \ \ \ \ \ \ \ \ \ $%
\begin{array}
[c]{c}
\end{array}
$

\textbf{Proof of Theorem 2.1. Necessity.} Let A and B be the losing score
lists of an [h-k]-bipartite hypertournament H(U, V). For any p and q with h
$\leq$\ p $\leq$\ m and k $\leq$\ q $\leq$\ n, let U$_{1}$ = \{u$_{1}$,
u$_{2}$,..., u$_{p}$\} and V$_{1}$ = \{v$_{1}$, v$_{2}$,..., v$_{q}$\} be the
set of vertices such that d$^{-}$(u$_{i}$) = a$_{i}$ for each 1 $\leq$\ i
$\leq$\ p, and d$^{-}$(v$_{j}$) = b$_{j}$ for each 1 $\leq$\ j $\leq$\ q. Let
H$_{1}$ be the [h-k]-bipartite subhypertournament formed by U$_{1}$ and
V$_{1}$. Then $\ \smallskip\ \ \ \ \ \ \ \ \ \ \ \ \ \ \ \sum_{i=1}^{p}%
a_{i}+\sum_{j=1}^{q}b_{j}\geq\sum_{i=1}^{p}d_{H_{1}}^{-}(u_{i})+\sum_{j=1}%
^{q}d_{H_{1}}^{-}(v_{j})=\left(
\begin{array}
[c]{c}%
p\\
h
\end{array}
\right)  \left(
\begin{array}
[c]{c}%
q\\
k
\end{array}
\right)  .$

\textbf{Sufficiency.} We induct on m and keep n fixed. For m = h, the result
is obviously true. Therefore, let m
$>$%
h, and similarly n
$>$%
k.

Now, \ \ \ \ $a_{m}=\sum_{i=1}^{m}a_{i}+\sum_{j=1}^{n}b_{j}-\left(  \sum
_{i=1}^{m-1}a_{i}+\sum_{j=1}^{n}b_{j}\right)  $

$\smallskip\ \ \ \ \ \ \ \ \ \ \ \ \ \ \ \ \ \ \leq\left(
\begin{array}
[c]{c}%
m\\
h
\end{array}
\right)  \left(
\begin{array}
[c]{c}%
n\\
k
\end{array}
\right)  -\left(
\begin{array}
[c]{c}%
m-1\\
h
\end{array}
\right)  \left(
\begin{array}
[c]{c}%
n\\
k
\end{array}
\right)  \smallskip$

$\smallskip\ \ \ \ \ \ \ \ \ \ \ \ \ \ \ \ \ \ =$\ $\left[  \left(
\begin{array}
[c]{c}%
m\\
h
\end{array}
\right)  -\left(
\begin{array}
[c]{c}%
m-1\\
h
\end{array}
\right)  \right]  \left(
\begin{array}
[c]{c}%
n\\
k
\end{array}
\right)  $

$\smallskip\ \ \ \ \ \ \ \ \ \ \ \ \ \ \ \ \ \ =$\ $\left(
\begin{array}
[c]{c}%
m-1\\
h-1
\end{array}
\right)  \left(
\begin{array}
[c]{c}%
n\\
k
\end{array}
\right)  $.

We consider the following two cases.

\textbf{Case 1. }\ $a_{m}=\left(
\begin{array}
[c]{c}%
m-1\\
h-1
\end{array}
\right)  \left(
\begin{array}
[c]{c}%
n\\
k
\end{array}
\right)  .$

So, \ \ \ $\sum_{i=1}^{m-1}a_{i}+\sum_{j=1}^{n}b_{j}=\sum_{i=1}^{m}a_{i}%
+\sum_{j=1}^{n}b_{j}-a_{m}$

$\smallskip
\ \ \ \ \ \ \ \ \ \ \ \ \ \ \ \ \ \ \ \ \ \ \ \ \ \ \ \ \ \ =\left(
\begin{array}
[c]{c}%
m\\
h
\end{array}
\right)  \left(
\begin{array}
[c]{c}%
n\\
k
\end{array}
\right)  -\left(
\begin{array}
[c]{c}%
m-1\\
h-1
\end{array}
\right)  \left(
\begin{array}
[c]{c}%
n\\
k
\end{array}
\right)  $

$\smallskip\ \ \ \ \ \ \ \ \ \ \ \ \ \ \ \ \ \ \ \ \ \ \ \ \ \ \ \ \ \ =$%
\ $\left[  \left(
\begin{array}
[c]{c}%
m\\
h
\end{array}
\right)  -\left(
\begin{array}
[c]{c}%
m-1\\
h-1
\end{array}
\right)  \right]  \left(
\begin{array}
[c]{c}%
n\\
k
\end{array}
\right)  =$\ $\left(
\begin{array}
[c]{c}%
m-1\\
h
\end{array}
\right)  \left(
\begin{array}
[c]{c}%
n\\
k
\end{array}
\right)  $.

By induction hypothesis [a$_{1}$, a$_{2}$,..., a$_{m-1}$] and B are losing
score lists of an [h-k]-bipartite hypertournament H$^{/}$(U$^{/}$, V) of order
m-1 x n. Construct an [h-k]-bipartite hypertournament H\ of order m x n as
follows. In H$^{/}$, let U$^{/}$ = \{u$_{1}$, u$_{2}$,..., u$_{m-1}$\} and V =
\{v$_{1}$, v$_{2}$,..., v$_{n}$\}. Adding a new vertex u$_{m}$, for each
(h+k)-tuple containing u$_{m}$, arrange u$_{m}$ on the last entry. Denote
E$_{1}$ to be the set of all these $\left(
\begin{array}
[c]{c}%
m-1\\
h-1
\end{array}
\right)  \left(
\begin{array}
[c]{c}%
n\\
k
\end{array}
\right)  $ (h+k)-tuple. Let E(H) = E(H$^{/}$)$\cup$E$_{1}$. Clearly, A and B
are losing score lists of H.

\textbf{Case 2. }\ $a_{m}<\left(
\begin{array}
[c]{c}%
m-1\\
h-1
\end{array}
\right)  \left(
\begin{array}
[c]{c}%
n\\
k
\end{array}
\right)  .$

Applying Lemma 2.3 repeatedly on A and keeping B fixed until we get a new
non-decreasing list A$^{/}=\ $[a$_{1}^{/}$, a$_{2}^{/}$,..., a$_{m}^{/}$] in
which now \ $a_{m}^{/}=\left(
\begin{array}
[c]{c}%
m-1\\
h-1
\end{array}
\right)  \left(
\begin{array}
[c]{c}%
n\\
k
\end{array}
\right)  .$ By Case 1, A$^{/}$ and B are the losing score lists of an
[h-k]-bipartite hypertournament. Now, apply Lemma 2.2 on A$^{/}$ and B
repeatedly until we obtain the initial pair of non-decreasing lists A and B.
Then by Lemma 2.2, A and B are the losing score lists of an [h-k]-bipartite
hypertournament. \ \ \ \ \ \ \ \ \ \ \ \ $%
\begin{array}
[c]{c}
\end{array}
\bigskip$

\textbf{Remark.} If h = 1, k =\ 1, we get the definition of scores in
bipartite tournaments and Theorem 2.1 gives

\begin{center}
$\sum_{i=1}^{p}a_{i}+\sum_{j=1}^{q}b_{j}\geq\left(
\begin{array}
[c]{c}%
p\\
1
\end{array}
\right)  \left(
\begin{array}
[c]{c}%
q\\
1
\end{array}
\right)  =pq,$
\end{center}

which is the characterization of score lists due to Beineke and Moon [2].\bigskip

\textbf{Proof of Theorem 2.2. }Let [c$_{1}$, c$_{2}$,..., c$_{m}$] and
[d$_{1}$, d$_{2}$,..., d$_{n}$] be score lists of an [h-k]-bipartite
hypertournament H(U, V), where U = \{u$_{1}$, u$_{2}$,..., u$_{m}$\} and V =
\{v$_{1}$, v$_{2}$,..., v$_{n}$\}with d$_{H}^{+}(u_{i})=c_{i}$ \ for i = 1,
2,..., m, and d$_{H}^{+}(v_{j})=d_{j}$ \ for j = 1, 2,..., n. Clearly,
d$^{+}(u_{i})+$d$^{-}(u_{i})=\left(
\begin{array}
[c]{c}%
m-1\\
h-1
\end{array}
\right)  \left(
\begin{array}
[c]{c}%
n\\
k
\end{array}
\right)  $ and d$^{+}(v_{j})+$d$^{-}(v_{j})=\left(
\begin{array}
[c]{c}%
m\\
h
\end{array}
\right)  \left(
\begin{array}
[c]{c}%
n-1\\
k
\end{array}
\right)  .$

Let a$_{m+1-i}=\ $d$^{-}(u_{i})$ and b$_{n+1-j}=$ d$^{-}(v_{j}).$

Then [a$_{1}$, a$_{2}$,..., a$_{m}$] and [b$_{1}$, b$_{2}$,..., b$_{n}$] are
the losing score lists of H. Conversely, if [a$_{1}$, a$_{2}$,..., a$_{m}$]
and [b$_{1}$, b$_{2}$,..., b$_{n}$] are the losing score lists of H, then
[c$_{1}$, c$_{2}$,..., c$_{m}$] and [d$_{1}$, d$_{2}$,..., d$_{n}$] are the
score lists of H. Hence it is sufficient to show that conditions (1) and (2)
are equivalent provided

\begin{center}
c$_{i}+$ a$_{m+1-i}=\left(
\begin{array}
[c]{c}%
m-1\\
h-1
\end{array}
\right)  \left(
\begin{array}
[c]{c}%
n\\
k
\end{array}
\right)  $
\end{center}

\ \ \ \ and \ \ \ \ \ \ \ \ \ \ \ \ \ \ \ \ \ \ \ \ d$_{j}+$ b$_{n+1-j}%
=\left(
\begin{array}
[c]{c}%
m\\
h
\end{array}
\right)  \left(
\begin{array}
[c]{c}%
n-1\\
k
\end{array}
\right)  .$

First, assume (2) holds. Then

$\smallskip\sum_{i=1}^{p}a_{i}+\sum_{j=1}^{q}b_{j}$

$=\sum_{i=1}^{p}\left\{  \left(
\begin{array}
[c]{c}%
m-1\\
h-1
\end{array}
\right)  \left(
\begin{array}
[c]{c}%
n\\
k
\end{array}
\right)  -c_{m+1-i}\right\}  +\sum_{j=1}^{q}\left\{  \left(
\begin{array}
[c]{c}%
m\\
h
\end{array}
\right)  \left(
\begin{array}
[c]{c}%
n-1\\
k-1
\end{array}
\right)  -d_{n+1-j}\right\}  $

$=p\left(
\begin{array}
[c]{c}%
m-1\\
h-1
\end{array}
\right)  \left(
\begin{array}
[c]{c}%
n\\
k
\end{array}
\right)  +q\left(
\begin{array}
[c]{c}%
m\\
h
\end{array}
\right)  \left(
\begin{array}
[c]{c}%
n-1\\
k-1
\end{array}
\right)  -$

$\smallskip\ \ \ \ \ \ \ \ \ \ \ \ \ \left[  \sum_{i=1}^{m}c_{i}+\sum
_{j=1}^{n}d_{j}-\sum_{i=1}^{m-p}c_{i}-\sum_{j=1}^{n-q}d_{j}\right]  $

\ $\smallskip\geq p\left(
\begin{array}
[c]{c}%
m-1\\
h-1
\end{array}
\right)  \left(
\begin{array}
[c]{c}%
n\\
k
\end{array}
\right)  +q\left(
\begin{array}
[c]{c}%
m\\
h
\end{array}
\right)  \left(
\begin{array}
[c]{c}%
n-1\\
k-1
\end{array}
\right)  -(h+k-1)\left(
\begin{array}
[c]{c}%
m\\
h
\end{array}
\right)  \left(
\begin{array}
[c]{c}%
n\\
k
\end{array}
\right)  $

$\smallskip\ \ \ \ \ \ \ \ \ \ \ \ +(m-p)\left(
\begin{array}
[c]{c}%
m-1\\
h-1
\end{array}
\right)  \left(
\begin{array}
[c]{c}%
n\\
k
\end{array}
\right)  +(n-q)\left(
\begin{array}
[c]{c}%
m\\
h
\end{array}
\right)  \left(
\begin{array}
[c]{c}%
n-1\\
k-1
\end{array}
\right)  +$

$\smallskip\ \ \ \ \ \ \ \ \ \ \ \ \ \ \ \ \ \left(
\begin{array}
[c]{c}%
m-(m-p)\\
h
\end{array}
\right)  \left(
\begin{array}
[c]{c}%
n-(n-q)\\
k
\end{array}
\right)  -\left(
\begin{array}
[c]{c}%
m\\
h
\end{array}
\right)  \left(
\begin{array}
[c]{c}%
n\\
k
\end{array}
\right)  $

$\smallskip\ \ \ =\left(
\begin{array}
[c]{c}%
p\\
h
\end{array}
\right)  \left(
\begin{array}
[c]{c}%
q\\
k
\end{array}
\right)  ,$

with equality when p = m and q = n. Thus, (1) holds.

Now, when (1) holds, using a similar argument as above, we can prove that (2) holds.

This completes the proof of the Theorem. \ \ \ \ \ \ \ \ $%
\begin{array}
[c]{c}
\end{array}
\bigskip$

\textbf{Corollary 2.1.} Given non-negative integers m, n, h and k with
$\ m\geq h>1\ $ and $\ n\geq k>1,$the non-decreasing sequences A = [a$_{i}%
$]$_{1}^{m}$ and B = [b$_{j}$]$_{1}^{n}$ of non-negative integers are the
losing score lists of an [h-k]-bipartite hypertournament if and only if for
each p and q,

$\ \ \ \ \ \ \ \sum_{i=1}^{p}a_{i}+\sum_{j=1}^{q}b_{j}\leq\left(
\begin{array}
[c]{c}%
m\\
h
\end{array}
\right)  \left(
\begin{array}
[c]{c}%
n\\
k
\end{array}
\right)  -\left(
\begin{array}
[c]{c}%
m-p\\
h
\end{array}
\right)  \left(
\begin{array}
[c]{c}%
n-q\\
k
\end{array}
\right)  ,$

\textbf{Proof.} This follows from Theorem 2.1.\bigskip\bigskip

\bigskip

\begin{center}
\bigskip\textbf{References}
\end{center}

[1] \ C. M. Bang and H. Sharp Jr., Score vectors of tournaments, J.
Combin.Theory Ser. B 26 (1) (1979) 81-84.

[2] \ L. W. Beineke and J. W. Moon, On bipartite tournaments and scores, Proc.
Fourth International Graph Theory Conference, Kalamazoo (1980\ 55-71.

[3] \ C. Berge, Graphs and hypergraphs, translated from French by E. Minieka,
North-Holland Mathematical Library\ \ 6, North-Holland Publishing Co.,
Amsterdam, London, (1973).

[4] \ Y. Koh and S. Ree, Score sequences of hypertournament matrices, J.Korea
Soc. Math. Educ. Ser. B: Pure and\ Appl. Math. 8 (2) (2001) 185-191.

[5] \ Y. Koh and S. Ree, On k-hypertournament matrices, Linear Algebra and its
Applications 373 (2003) 183-195.

[6] \ H. G. Landau, On dominance relations and the structure of animal
societies. III. The condition for a score structure, Bull. Math. Biophys.\ 15
( 1953 ) 143-148.

[7] \ C. Wang and G. Zhou, Note on the degree sequences of k-hypertournaments,
Discrete Mathematics, Preprint.

[8] \ G. Zhou, T. Yao and K. Zhang, On score sequences of k-hypertournaments,
European J. Combin. 21 (8)\ (2000) 993-1000.
\end{document}